\documentclass[11pt,A4paper]{article}
\usepackage{amsmath,amsfonts,mathrsfs,epsfig}

\newtheorem{thm}{Theorem}[section]

\newtheorem{lem}[thm]{Lemma}
\newtheorem{cor}[thm]{Corollary}

\newcommand{\maj}{\mathrm{maj}}
\newcommand{\fmaj}{\mathrm{fmaj}}
\newcommand{\Drg}{\mathscr{D}}

\def\pf{\noindent{\it Proof.} }
\def\qed{\nopagebreak\hfill{\rule{4pt}{7pt}}\medbreak}

\makeatletter \@addtoreset{equation}{section} \makeatother

\begin{document}

\begin{center}
{\Large\bf The Limiting Distributions of the Coefficients of

 the
$q$-Derangement Numbers}
\end{center}

\vskip 2mm

\centerline{William Y.C. Chen$^1$ and David G.L. Wang$^2$}

\begin{center}
Center for Combinatorics, LPMC-TJKLC\\
Nankai University, Tianjin 300071, P.R. China

\vskip 2mm%

$^1$chen@nankai.edu.cn, $^2$wgl@cfc.nankai.edu.cn
\end{center}

\begin{abstract}
We show that  the distribution of the coefficients of the
$q$-derangement numbers is asymptotically normal. We also show that
this property holds for the $q$-derangement numbers of type $B$.
\end{abstract}

\noindent\textbf{Keywords:} $q$-derangement numbers, flag major
index, moment generating function, limiting distribution

\noindent\textbf{AMS Classification:} 05A15, 05A16, 05A30

\section{Introduction}

Let $\mathfrak{S}_n$ denote the symmetric group of permutations on
$[n]=\{1,2,\ldots,n\}$.  Let $\Drg_n$ denote the set of {\em
derangements}, i.e.,
\[
\Drg_n=\{\pi=\pi_1\pi_2\cdots\pi_n\in\mathfrak{S}_n\colon\pi_i\ne
i,\ i=1,2,\ldots,n\}.
\]
The major index of a permutation $\pi=\pi_1\pi_2\cdots \pi_n$ is
defined by
\[
\maj(\pi)=\sum_{\pi_i>\pi_{i+1}}i.
\]
The following formula was derived by Gessel and published in
\cite{Ges-Reu93}:
\begin{equation}\label{d_n(q)_A}
d_n(q)%
=\sum_{\pi\in\Drg_n}q^{\maj(\pi)}
=[n]_q!\sum_{k=0}^n\frac{(-1)^kq^{k\choose2}}{[k]_q!},
\end{equation}
where $[0]_q=[0]_q!=1$ and for $k\ge1$,
$[k]_q=1+q+q^2+\cdots+q^{k-1}$ and $[k]_q!=[k]_q[k-1]_q\cdots[1]_q$.
The coefficients of $d_n(q)$   are given in Table \ref{t1} for $n\le
6$. Combinatorial proofs of \eqref{d_n(q)_A} have been found by
Wachs \cite{Wac89}, and Chen and Xu \cite{Chen-Xu06}.

\begin{table}[h]
\begin{center}
\begin{tabular}{|r||c|c|c|c|c|c|c|c|c|c|c|c|c|c|c|}
 \hline
 $n\backslash k$&1&2&3&4&5&6&7&8&9&10&11&12&13&14&15\\
 \hline\hline
 2&1&&&&&&&&&&&&&&\\
 \hline
 3&1&1&&&&&&&&&&&&&\\
 \hline
 4&1&2&2&2&1&1&&&&&&&&&\\
 \hline
 5&1&3&5&7&8&8&6&4&2&&&&&&\\
 \hline
 6&1&4&9&16&24&32&37&38&35&28&20&12&6&2&1\\
 \hline
\end{tabular}
\caption{\label{t1} The $q$-derangement numbers of type $A$ for
$n\le6$.}
\end{center}
\end{table}

\begin{figure}[h]
\begin{center}
\includegraphics[width=8cm]{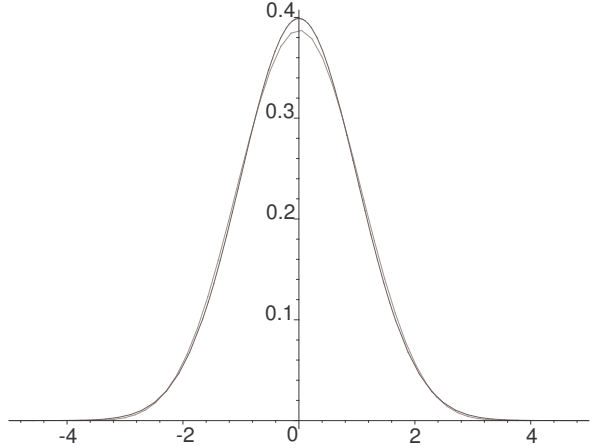}
\caption{\label{f1} The distribution of the coefficients of
$\Drg_{10}$ compared with the normal distribution.}
\end{center}
\end{figure}

\begin{figure}[h]
\begin{center}
\includegraphics[width=8cm]{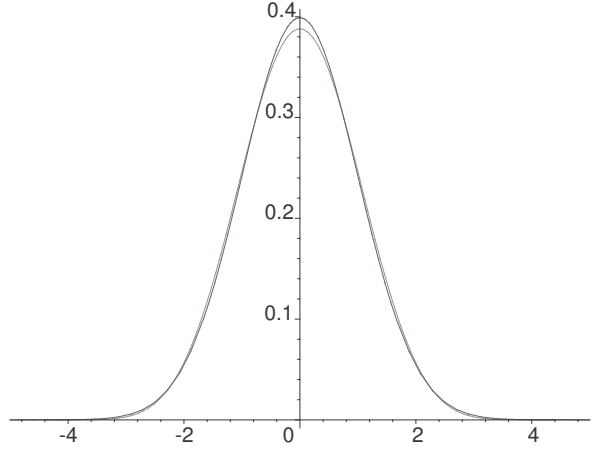}
\caption{\label{f2} The distribution of the coefficients of
$\Drg_{10}^B$ compared with the normal distribution.}
\end{center}
\end{figure}

In this paper, we will show that the limiting distribution of the
coefficients of $d_n(q)$, that is, the major index of a random
derangement, is normal, see Figure \ref{f1}. Moreover, we will show
that the limiting distribution of the $q$-derangement numbers of
type $B$ is also normal, see Figure \ref{f2}.

Write the set $\{\bar{1},\bar{2},\ldots,\bar{n}\}$ as $[\bar{n}]$.
Let $\mathfrak{S}_n^B$ denote the hyperoctahedral group of
permutations on $[n]\cup[\bar{n}]$, called signed permutations or
$B_n$-permutations, see Bj\"{o}rner and Brenti \cite{Bjo-Bre05}. Let
$\Drg_n^B$ denote the set of {\em $B_n$-derangements} on $[n]$,
namely,
\[
\Drg_n^B=\{\pi=\pi_1\pi_2\cdots\pi_n\in\mathfrak{S}_n^B\colon\pi_i\ne
i,\ i=1,2,\ldots,n\}.
\]
For example, $\Drg_1^B=\{\bar{1}\}$, $\Drg_2^B=\left\{
\bar{1}\bar{2}, 21, 2\bar{1}, \bar{2}1, \bar{2}\bar{1} \right\}$.
For $B_n$-permutations, Adin and Roichman \cite{Adi-Roi01}
introduced the notion of the {\em flag major index}, or the $\fmaj$
index for short, defined by
\[
\fmaj(\pi)=2\maj(\pi)+\mathrm{neg}(\pi),
\]
where $\maj(\pi)$ is the major index of $\pi$ with respect to the
following order on $[n]\cup [\bar{n}]$:
\[
\bar{n}<\cdots<\bar{2}<\bar{1}<1<2<\cdots<n,
\]
and $\mathrm{neg}(\pi)$ is the number of $\pi_i$'s in $[\bar{n}]$,
see also Adin, Brenti and Roichman \cite{Adi-Bre-Roi01}, and Chow
and Gessel \cite{Chow-Ges03}. For example, the flag major of the
$B_7$-permutation $35\bar{1}2\bar{6}\bar{7}4$ equals $2\times
11+3=25$.   Chow \cite{Chow06} derived the following formula for the
$q$-derangement numbers of type $B$:
\begin{equation}\label{d_n(q)_B}
d_n^B(q)%
=\sum_{\pi\in\Drg_n^B}q^{\fmaj(\pi)}%
=[2n]_q!!\sum_{k=0}^n\frac{(-1)^kq^{k(k-1)}}{[2k]_q!!},
\end{equation}
where $[2k]_q!!=[2k]_q[2k-2]_q\cdots[2]_q$. For $n\le 4$, the
coefficients of the polynomials $d_n^B(q)$ are given in Table
\ref{t2}.
\begin{table}[h]
\begin{center}
\begin{tabular}{|r||c|c|c|c|c|c|c|c|c|c|c|c|c|c|c|c|}
 \hline
 $n\backslash k$&1&2&3&4&5&6&7&8&9&10&11&12&13&14&15&16\\
 \hline\hline
 1&1&&&&&&&&&&&&&&&\\
 \hline
 2&1&2&1&1&&&&&&&&&&&&\\
 \hline
 3&1&3&4&5&5&4&4&2&1&&&&&&&\\
 \hline
 4&1&4&8&13&18&22&26&28&28&25&21&17&11&7&3&1\\
 \hline
\end{tabular}
\caption{\label{t2} The $q$-derangement numbers of type $B$ for
$n\le4$.}
\end{center}
\end{table}

Based on the formula (\ref{d_n(q)_B}), we will show that the
limiting distribution of the coefficients of $d_n^B(q)$ is normal.
 Figure \ref{f2} is an illustration of the distribution  for $n=10$.

\section{The Limiting Distribution of the Coefficients of $d_n(q)$}%
\label{Section_A}

The aim of this section is to show that the limiting distribution of
the coefficients of $d_n(q)$ is normal. We write
\[
f_{n,k}(q)%
=\begin{cases}
1, &\mbox{if } k=n;\\[5pt]
[n]_q[n-1]_q\cdots[k+1]_q, &\mbox{else}.
\end{cases}
\]
Then we can express $d_n(q)$ as
\begin{equation}\label{GF3}
d_n(q)=\sum_{k=0}^n(-1)^kq^{k\choose2}f_{n,k}(q).
\end{equation}
Let $D_n=|\Drg_n|$ be the number of derangements in $\Drg_n$. For
example, $D_1=0$, $D_2=1$, $D_3=2$, $D_4=9$, $D_5=44$.

We will adopt the common notation in asymptotic analysis. If $f(n)$
and $g(n)$ are two functions of $n$, then
\begin{itemize}
\item%
$f(n)=o(g(n))$ means that $\lim_{n\to\infty}f(n)/g(n)=0$;
\item%
$f(n)\sim g(n)$ means that $\lim_{n\to\infty}|f(n)|/|g(n)|=1$.
\end{itemize}

We now recall some basic facts about the derangement numbers $D_n$,
see, for example, Stanley \cite{Sta97}. For $n\ge3$,
\begin{align}
D_n%
&=n!\sum_{k=0}^n\frac{(-1)^k}{k!}\label{D_n_A1}\\
&=nD_{n-1}+(-1)^n\label{D_n_A2}\\[9pt]
&=(n-1)(D_{n-1}+D_{n-2})\label{D_n_A3}\\[5pt]
&=\left\lfloor\frac{n!}{e}+\frac{1}{2}\right\rfloor%
\sim\frac{n!}{e},\label{D_n_A4}
\end{align}
where the symbol $\lfloor x\rfloor$ denotes the largest integer not
exceeding $x$. From \eqref{D_n_A2} it immediately follows that
\begin{equation}\label{Ratio_D_A}
\frac{D_{n-1}}{D_n}%
=\frac{1}{n}-\frac{(-1)^n}{nD_n}%
\frac{1}{n}+o(1).
\end{equation}

While it is common to use $\maj(\pi)$ to denote the major index of a
permutation $\pi$, there does not seem to be any confusion if we
also use $\maj$ to denote the major index of a random derangement on
$[n]$. The probability generating function of $\maj$ is clearly
$d_n(x)/D_n$, whereas the moment generating function of $\maj$ is
given by
\begin{equation}\label{MGF_A}
M_n(x)=\frac{1}{D_n}d_n\left(e^x\right)%
=\frac{1}{D_n}[n]_{e^x}!%
\sum_{k=0}^n\frac{(-1)^k\,e^{x{k\choose2}}}{[k]_{e^x}!}.
\end{equation}
Let $E_n$, $V_n$ and $\sigma_n=V_n^{1/2}$ denote the expectation,
the variance and the standard deviation of $\maj$ respectively. Then
the probability generating function $\widetilde{d_n}(q)$ of the
normalized random variable $(\maj-E_n)/\sigma_n$ equals
\[
\widetilde{d_n}(q)%
=\sum_{\pi\in\Drg_n}q^{(\maj(\pi)-E_n)/\sigma_n}%
=q^{-E_n/\sigma_n}d_n(q^{1/\sigma_n}).
\]
Thus by the definition \eqref{MGF_A}, the moment generating function
of $(\maj-E_n)/\sigma_n$ equals
\begin{equation}\label{MGF_Standard_A}
\widetilde{M_n}(t)=\widetilde{d_n}(e^t)/D_n%
=\exp(-t\,E_n/\sigma_n)M_n(t/\sigma_n).
\end{equation}

\subsection{The expectation and variance}

We now compute the expectation and variance of the major index
$\maj$ of a random derangement on $[n]$.

\begin{thm}\label{EV_A}%
The expectation $E_n$ and variance $V_n$ of the random variable
$\maj$ given by
\begin{equation}\label{E_A}
E_n%
=\frac{1}{2}{n\choose2}\left(1+\frac{D_{n-2}}{D_n}\right)%
=\frac{n^2-n+1}{4}+\frac{(-1)^n(n-1)}{4D_n},
\end{equation}
and
\begin{equation}\label{V_A}
V_n=\frac{2n^3+3n^2-5n-16}{72}%
+\frac{9n^3-4n^2-46n+41}{144}\frac{(-1)^n}{D_n}%
-\left(\frac{n-1}{4D_n}\right)^2.
\end{equation}
\end{thm}

Here we give only a sketch of the proof, and detailed steps are
omitted.

\pf The generating function \eqref{d_n(q)_A} implies that
\begin{align}
E_n%
&=\frac{1}{D_n}\sum_{\pi\in\Drg_n}\maj(\pi)=\frac{d_n'(1)}{D_n},%
\label{E_A_tmp}\\[5pt]
V_n%
&=\frac{1}{D_n}\sum_{\pi\in\Drg_n}\maj^2(\pi)-E_n^2%
=\frac{d_n''(1)}{D_n}+E_n-E_n^2,\label{V_A_tmp}
\end{align}
where $d'_n(q)$ and $d''_n(q)$ are the first and second derivatives
of $d_n(q)$. From \eqref{GF3} and \eqref{D_n_A1}, we find
\[
d_n'(1)%
=\sum_{k=0}^n(-1)^k\left[{k\choose2}f(1)+f'(1)\right]%
=\frac{1}{2}{n\choose2}(D_n+D_{n-2}).
\]
So \eqref{E_A} follows from \eqref{E_A_tmp}. Differentiating
\eqref{GF3} twice yields
\begin{equation}\label{d''}
d_n''(1)%
=\sum_{k=0}^n(-1)^k\left[%
{k\choose2}\left({k\choose2}-1\right)f(1)+2{k\choose2}f'(1)+f''(1)%
\right].
\end{equation}
The following relations can be easily verified:
\[
\sum_{k=0}^n(-1)^k\left[%
{k\choose2}\left({k\choose2}-1\right)f(1)+2{k\choose2}f'(1)\right]%
=\frac{3}{2}{n\choose3}(n+1)D_{n-2},
\]
\begin{align*}
\sum_{k=0}^n(-1)^kf''(1)%
={n\choose3}\frac{9(n-3)D_{n-4}-32D_{n-3}-18(n+1)D_{n-2}+(9n+13)D_n}{24}.
\end{align*}
Now, using \eqref{d''} and \eqref{D_n_A3}, we deduce that
\[
d_n''(1)=\frac{1}{72}{n\choose2}%
[(n-2)(27n+32)D_{n-2}-(9n+5)D_{n-1}+(n-2)(9n+13)D_n].
\]
According to \eqref{V_A_tmp},
\begin{align*}
V_n&=\frac{n}{144}\left[%
-(n-1)\left(27n^2-13n-23\right)\frac{D_{n-1}}{D_n}%
+\left(9n^3+13n^2-7n-38\right)\right]\\[5pt]%
&\quad-\left(\frac{n^2}{4}-\frac{n(n-1)}{4}\frac{D_{n-1}}{D_n}\right)^2.
\end{align*}
In view of \eqref{D_n_A3} and \eqref{Ratio_D_A}, we obtain
\eqref{V_A}.\qed

We note that the formula \eqref{E_A} for the expectation of the
major index can also be derived by a combinatorial argument, the
details are omitted. Based on the estimates \eqref{D_n_A4} and
\eqref{Ratio_D_A}, we derive the following approximations.

\begin{cor}\label{EV_Est_A}
We have the following asymptotic estimates:
\[
E_n=\frac{n^2}{4}-\frac{n}{4}+\frac{1}{4}+o(1),\quad\quad
V_n=\frac{n^3}{36}+\frac{n^2}{24}-\frac{5n}{72}-\frac{2}{9}+o(1).
\]
\end{cor}

\subsection{The limiting distribution}
It is well-known that the moment generating function of a random
variable determines its distribution by Curtiss's theorem (see
Curtiss \cite{Cur42} or Sachkov \cite{Sac97}). In particular, if the
moment generating function $M_n(x)$ of a random variable $\xi_n$ has
the limit
\[
\lim_{n\to\infty}M_n(x)=e^{x^2/2},
\]
then $\xi_n$ has as an asymptotically standard normal distribution
as $n$ trending to infinity.

We will need Tannery's theorem (see Tannery \cite{Tan1904}) which is
essential in the proofs of Lemma \ref{lem_Tannery_A} and Lemma
\ref{lem_Tannery_B}.

\begin{thm}[Tannery's theorem]
Let $\{v_k(n)\}_{k\ge0}$ be an infinite series satisfying the
following two conditions.
\begin{itemize}
\item%
For any fixed $k$, there holds $\lim_{n\to\infty}v_k(n)=w_k$.
\item%
For any non-negative integer $k$, $|v_k(n)|\le M_k$, where $M_k$
independent of $n$ and the series $\sum_{k\ge0}M_k$ is convergent.
\end{itemize}
Then
\[
\lim_{n\to\infty}\sum_{k=0}^{m(n)}v_k(n)=\sum_{k=0}^\infty w_k,
\]
where $m(n)$ is an increasing integer-valued function which trends
steadily to infinity as $n$ does.
\end{thm}

\begin{lem}\label{lem_Tannery_A}
For any $|x|\le1$ and bounded $|t|\le M$, we have
\begin{equation}\label{eq9}
\lim_{n\to\infty}\sum_{k=0}^n\frac{x^k}{[k]_{e^{-t/\sigma_n}}!}=e^x.
\end{equation}
\end{lem}

\pf We apply Tannery's theorem and set
\[
v_k(n)=\frac{x^k}{[k]_{e^{-t/\sigma_n}}!},
\]
and $m(n)=n$. Then for any fixed $k$, by Corollary \ref{EV_Est_A},
it is clear that
\[
w_k=\lim_{n\to\infty}v_k(n)=\frac{x^k}{k!}.
\]
Note that the right hand side of \eqref{eq9} can be expressed as
\[
e^x=\sum_{k=0}^\infty\frac{x^k}{k!}=\sum_{k=0}^\infty w_k.
\]
By virtue of Tannery's theorem, to prove \eqref{eq9} it suffices to
find an upper bound $M_k$ for
\[
|v_k(n)|=\frac{|x^k|}{[k]_{e^{-t/\sigma_n}}!},
\]
such that $M_k$ is independent of $n$ and $\sum_{k=0}^{\infty}M_k$
converges. We claim that there exists a constant $c\in(0,1]$ such
that $M_k=(1+c)^{1-k}$ is the desired upper bound and this bound
clearly implies the convergence of $\sum_{k=0}^\infty M_k=1/c+2+c$.

For $t\le0$, we have $e^{-t/\sigma_n}\ge1$ and thus
\[
\frac{|x^k|}{[k]_{e^{-t/\sigma_n}}!}%
\le\frac{|x^k|}{k!}\le\frac{1}{k!}\le\frac{1}{2^{k-1}}\le M_k.
\]
For $t\ge0$, Corollary \ref{EV_Est_A} implies that $\sigma_n$ has a
positive lower bound as $n$ runs over all positive integers and so
does $e^{-t/\sigma_n}$. Suppose that $e^{-t/\sigma_n}\ge
c_t\in(0,1]$. Since the function $e^{-t/\sigma_n}$ is continuous in
$t$ and $t$ is bounded, there exists a constant $c\in(0,1]$
independent of $t$ so that $e^{-t/\sigma_n}\ge c$ for all $|t|\le
M$. Hence for any $k\ge1$,
\begin{align*}
\frac{|x^k|}{[k]_{e^{-t/\sigma_n}}!}%
&=\prod_{j=1}^k\frac{|x|}{1+e^{-t/\sigma_n}+%
\cdots+e^{-(j-1)t/\sigma_n}}\\[5pt]
&\le\prod_{j=1}^k\frac{1}{1+c+\cdots+c^{j-1}}\\[5pt]
&\le\prod_{j=2}^k\frac{1}{1+c}=M_k.
\end{align*}
This completes the proof. \qed

In the computation of the moment generating function of $\maj$, we
will need the Bernoulli numbers $B_k$ which have the following
generating function,
\begin{equation}\label{Bernoulli_GF}
\frac{x}{e^x-1}=\sum_{k=0}^{\infty}B_k\frac{x^k}{k!}.
\end{equation}
The first few Bernoulli numbers are
\[
B_0=1,\ B_1=-1/2,\ B_2=1/6,\ B_3=0,\ B_4=-1/30.
\]
Moreover, $B_{2i+1}=0$ for any $i\ge1$. Alzer \cite{Alz00}
establishes sharp bounds for $|B_{2n}|$ leading to the following
asymptotic formula (see also \cite[pp. 805]{Abr-Ste64}) which will
be needed in the proof of Lemma \ref{lem_i>=2_A}:
\begin{equation}\label{Bernoulli}
|B_{2n}|\sim\frac{2\cdotp(2n)!}{(2\pi)^{2n}}.
\end{equation}

\begin{lem}\label{lem_i>=2_A}%
For any bounded $|t|<M$, we have
\begin{equation} \label{eq1}
\lim_{n\to\infty}\sum_{i=2}^{\infty}%
\frac{B_{2i}\,t^{2i}}{(2i)\,(2i)!\,\sigma_n^{2i}}%
\sum_{j=1}^n\left(j^{2i}-1\right)=0,
\end{equation}
where $B_{2i}$ are the Bernoulli numbers.
\end{lem}

\pf Let $\alpha$, $\beta$ and $\gamma$ be three constants such that
$\alpha>1$, $\beta>36$, and $0<\gamma<1/2$. Let $N$ be a fixed
integer satisfying the following three conditions:
\begin{itemize}
\item%
$n+1<\alpha n$ for any $n>N$;
\item%
$\sigma_n^2-n^3/\beta>0$ for any $n>N$;
\item%
$2\pi N^{\gamma/2}>M\alpha\sqrt{\beta}$.
\end{itemize}
The existence of such $N$ is obvious. Let $i\ge2$ and $n>N$. From
the inequalities
\[
\sum_{j=1}^{n}\left(j^{2i}-1\right)%
<\int_1^{n+1}\left(t^{2i}-1\right)dt%
=\frac{(n+1)^{2i+1}-1}{2i+1}-n%
<\frac{(n+1)^{2i+1}}{5}%
<\frac{(\alpha n)^{2i+1}}{5}
\]
and the assumption $\sigma_n^2>n^3/\beta$, we deduce that
\[
\frac{1}{\sigma_n^{2i}}\sum_{j=1}^{n}\left(j^{2i}-1\right)%
<\frac{\beta^i}{n^{3i}}%
\frac{(\alpha n)^{2i+1}}{5}%
=\frac{\alpha}{5}%
\frac{\left(\alpha\sqrt{\beta}\right)^{2i}}{n^{i-1}}.
\]
In light of the inequality
\[
\frac{1}{n^{i-1}}%
=\frac{1}{n^{\gamma i}}\cdot\frac{1}{n^{(1-\gamma)i-1}}%
<\frac{1}{N^{\gamma i}}\cdot\frac{1}{n^{1-2\gamma}},
\]
we see that
\begin{equation}\label{eq2}
\lim_{n\to\infty}\sum_{i=2}^{\infty}%
\frac{B_{2i}\,t^{2i}}{(2i)\,(2i)!\,\sigma_n^{2i}}%
\sum_{j=1}^n\left(j^{2i}-1\right)%
\le\frac{\alpha}{5}\lim_{n\to\infty}%
\left(\sum_{i=2}^{\infty}\frac{|B_{2i}|}{(2i)\,(2i)!}%
\frac{(\alpha\sqrt{\beta})^{2i}}{N^{\gamma i}}t^{2i}\right)%
n^{2\gamma-1}.
\end{equation}
By the asymptotic estimate \eqref{Bernoulli} for Bernoulli numbers,
we see that the radius of convergence (see, for example, Howie
\cite{How03}) of the series on the right hand of \eqref{eq2} equals
\[
\lim_{i\to\infty}%
\left(\frac{|B_{2i}|}{(2i)\,(2i)!}%
\frac{\left(\alpha\sqrt{\beta}\right)^{2i}}%
{N^{\gamma i}}\right)^{-\frac{1}{2i}}%
=\frac{2\pi N^{\gamma/2}}{\alpha\sqrt{\beta}}>M.
\]
Since $\lim_{n\to\infty}n^{2\gamma-1}=0$, we conclude that the
series in \eqref{eq1} is absolutely convergent to zero for
$|t|<M$.\qed

The following lemma gives an expression of the moment generating
function of the random variable $\maj$ in term of the Bernoulli
numbers. This lemma will be needed in the proof of Theorem
\ref{Main_A}.

\begin{lem}\label{lem_MGF_A}
The moment generating function of $\maj$ equals
\[
M_n(x)=\frac{n!}{D_n}\exp\left(\frac{n(n-1)x}{4}+%
\sum_{i=1}^\infty\frac{B_{2i}\,x^{2i}}{(2i)\,(2i)!}%
\sum_{j=1}^n\left(j^{2i}-1\right)\right)%
\sum_{k=0}^n\frac{(-1)^k}{[k]_{e^{-x}}!}.
\]
\end{lem}

\pf By the formula \eqref{MGF_A}, we need to express $[n]_{e^x}!$
and $e^{x{k\choose2}}\big/[k]_{e^x}!$ in terms of Bernoulli numbers.
It is known that, see, for example, Mcintosh \cite{Mci99},
\[
1-e^{-x}%
=x\cdotp\exp\left(\sum_{k=1}^\infty\frac{B_n\,x^k}{k\cdotp
k!}\right).
\]
Thus for any $j\ge1$,
\begin{align*}
1-e^{xj}%
&=-xj\cdotp%
\exp\left(\sum_{i=1}^\infty\frac{B_i(-xj)^i}{i\cdotp i!}\right)%
=-xj\cdotp%
\exp\left(\frac{xj}{2}+\sum_{i=1}^\infty%
\frac{B_{2i}(xj)^{2i}}{(2i)\,(2i)!}\right),\\[5pt]
[j]_{e^x}%
&=\frac{1-e^{xj}}{1-e^x}%
=j\cdotp\exp\left(\frac{x(j-1)}{2}+%
\sum_{i=1}^\infty\frac{B_{2i}\,x^{2i}\left(j^{2i}-1\right)}%
{(2i)\,(2i)!}\right).
\end{align*}
Therefore,
\begin{align}
[n]_{e^x}!%
&=\prod_{j=1}^n[j]_{e^x}%
=\prod_{j=1}^nj\cdotp\exp\left(\frac{x(j-1)}{2}%
+\sum_{i=1}^\infty\frac{B_{2i}\,x^{2i}\left(j^{2i}-1\right)}%
{(2i)\,(2i)!}\right)\notag\\[5pt]
&=n!\cdotp\exp\left(\frac{n(n-1)x}{4}+%
\sum_{i=1}^\infty\frac{B_{2i}\,x^{2i}}{(2i)\,(2i)!}%
\sum_{j=1}^n\left(j^{2i}-1\right)\right).\label{eq3}
\end{align}
Observe that
\begin{equation}\label{eq4}
\frac{e^{x{k\choose2}}}{[k]_{e^x}!}%
=e^{x{k\choose2}}%
\left(\prod_{j=1}^k\frac{1-e^{jx}}{1-e^x}\right)^{-1}
=\prod_{j=1}^k\frac{1-e^x}{1-e^{jx}}\frac{e^{jx}}{e^x}
=\frac{1}{[k]_{e^{-x}}!}.
\end{equation}
Substituting \eqref{eq3} and \eqref{eq4} into \eqref{MGF_A}, we
obtain the desired expression. \qed

\begin{thm}\label{Main_A}
Let $\maj$ be the major index of a random derangement on $[n]$. Then
the distribution of the random variable
\[
\xi_n={\maj-E_n \over \sigma_n}
\]
converges to the standard normal distribution as $n\to\infty$.
\end{thm}

\pf By Curtiss's theorem and \eqref{MGF_Standard_A}, the normality
of the distribution of the standardized random variable $\xi_n$ can
be justified by the following relation
\[
\lim_{n\to\infty}e^{-t\,E_n/\sigma_n}M_n(t/\sigma_n)=e^{t^2/2}.
\]
By virtue of Lemma \ref{lem_MGF_A}, the above relation can be
restated as
\[
\lim_{n\to\infty}
\frac{n!}{D_n}\exp\left(-\frac{tE_n}{\sigma_n}%
+\frac{n(n-1)t}{4\sigma_n}%
+\sum_{i=1}^\infty\frac{B_{2i}\,t^{2i}}{(2i)\,(2i)!\sigma_n^{2i}}%
\sum_{j=1}^n\left(j^{2i}-1\right)\right)%
\sum_{k=0}^n\frac{(-1)^k}{[k]_{e^{-t/\sigma_n}}!}%
=e^{t^2/2}.
\]
First of all, the estimate \eqref{D_n_A4} implies that
\begin{equation}\label{eq1_A}
\lim_{n\to\infty}n!/D_n=e.
\end{equation}
By Corollary \ref{EV_Est_A}, for bounded $t$ we have
\begin{equation}\label{eq2_A}
\lim_{n\to\infty}%
\left(\frac{n(n-1)t}{4\sigma_n}-\frac{tE_n}{\sigma_n}\right)%
=\lim_{n\to\infty}%
\frac{t}{\sigma_n}\left(\frac{n(n-1)}{4}-E_n\right)%
=0,
\end{equation}
It is easily checked that
\[
\lim_{n\to \infty} {1\over \sigma_n^2}
\sum_{j=1}^n\left(j^2-1\right) = 12.
\]
In view of Lemma \ref{lem_i>=2_A} and the fact that $B_2=1/6$, we
have
\begin{equation}\label{eq3_A}
\lim_{n\to\infty}%
\sum_{i=1}^\infty\frac{B_{2i}\,t^{2i}}{(2i)\,(2i)!\sigma_n^{2i}}%
\sum_{j=1}^n\left(j^{2i}-1\right)%
=\lim_{n\to\infty}%
\frac{B_2\,t^2}{2\cdot2!\,\sigma_n^{2}}\sum_{j=1}^n\left(j^2-1\right)%
=\frac{t^2}{2}.
\end{equation}
Finally, taking $x=-1$ in Lemma \ref{lem_Tannery_A}, we get
\begin{equation}\label{eq4_A}
\lim_{n\to\infty}\sum_{k=0}^n\frac{(-1)^k}{[k]_{e^{-t/\sigma_n}}!}%
=e^{-1}.
\end{equation}
Combining \eqref{eq1_A}, \eqref{eq2_A}, \eqref{eq3_A} and
\eqref{eq4_A}, we complete the proof. \qed

\section{The Limiting Distribution of the Coefficients of $d_n^B(q)$}%
\label{Section_B}

In this section, we show that the limiting distribution of the
$q$-derangement numbers is normal. Let $D_n^B$ be the number of
$B_n$-derangements on $[n]$. The first few values of $D_n^B$ are
\[
D_1^B=1,\ D_2^B=5,\ D_3^B=29,\ D_4^B=233,\ D_5^B=2329,\ D_6^B=27949.
\]
For $n\ge3$, we have
\begin{align}
D_n^B%
&=d_n^B(1)=(2n)!!\sum_{k=0}^n\frac{(-1)^k}{(2k)!!}\label{D_n_B1}\\
&=2nD_{n-1}^B+(-1)^n\label{D_n_B2}\\[10pt]
&=(2n-1)D_{n-1}^B+(2n-2)D_{n-2}^B\label{D_n_B3}\\[5pt]
&=\left\lfloor\frac{(2n)!!}{\sqrt{e}}+\frac{1}{2}\right\rfloor%
\sim\frac{(2n)!!}{\sqrt{e}}. \label{D_n_B4}
\end{align}
For completeness, we present a proof for \eqref{D_n_B4}:
\[
D_n^B%
=(2n)!!\sum_{k=0}^n\left(-\frac{1}{2}\right)^k\frac{1}{k!}%
=(2n)!!\left(e^{-1/2}-\sum_{k=n+1}^{\infty}\frac{(-1)^k}{(2k)!!}\right).
\]
It is easy to see that the absolute value of the remainder
\[
r_n=\sum_{k=n+1}^{\infty}\frac{(-1)^k}{(2k)!!}
\]
is not greater than the absolute value of the $(n+1)$-st term of the
alternating series, i.e., $1/(2n+2)!!$. This yields
\[
D_n^B%
=(2n)!!/\sqrt{e}-y_n,
\]
where
\[
|y_n|=|(2n)!!r_n|\le(2n+2)^{-1}\le1/4.
\]
Since $D_n^B$ is an integer, \eqref{D_n_B4} is verified.

From \eqref{D_n_B2} it follows that
\begin{equation}\label{Ratio_D_B}
\frac{D_{n-1}^B}{D_n^B}%
=\frac{1}{2n}-\frac{(-1)^n}{2nD_n^B}%
\frac{1}{2n}+o(1).
\end{equation}

Let $E_n^B$, $V_n^B$ and $\sigma_n^B=\left(V_n^B\right)^{1/2}$
denote the expectation, the variance and the standard deviation of
$\fmaj$ respectively. We also use $\fmaj$ to denote the fmaj index
of a random $B_n$-derangements on $[n]$. The probability generating
function of $\fmaj$ is
\[
\widetilde{d_n^B}(q)%
=q^{-E_n^B/\sigma_n^B}d_n^B\left(q^{1/\sigma_n^B}\right).
\]
The moment generating function of $\fmaj$ is given by
\begin{equation}\label{MGF_B}
M_n^B(x)%
=\frac{1}{D_n^B}d_n^B\left(e^x\right)%
=\sum_{k=0}^\infty\left(\frac{1}{D_n^B}%
\sum_{\pi\in\Drg_n^B}\fmaj(\pi)^k\right).
\end{equation}
The normalized random variable $(\fmaj-E_n^B)/\sigma_n^B$ equals
\begin{equation}\label{MGF_Standard_B}
\widetilde{M_n^B}(t)%
=\exp\left(-t\,E_n^B/\sigma_n^B\right)M_n^B\left(t/\sigma_n^B\right).
\end{equation}

\subsection{The expectation and variance}

Let $[x^i]f(x)$ to denote the coefficient of $x^i$ in the expansion
of $f(x)$. Then the expectation and variance of $\fmaj$ can be
expressed in terms of the moment generating function $M_n^B(x)$:
\begin{align}
E_n^B&=\frac{1}{D_n^B}\sum_{\pi\in\Drg_n^B}\fmaj=[x]M_n^B(x).%
\label{E_B_tmp}\\[5pt]
V_n^B%
&=\left(\frac{1}{D_n^B}%
\sum_{\pi\in\Drg_n^B}\fmaj^2\right)-\left(E_n^B\right)^2%
=2\left[x^2\right]M_n^B(x)-\left(E_n^B\right)^2.\label{V_B_tmp}
\end{align}

Let $\left\langle x^2 \right\rangle f(x)$ to denote the truncated
sum of $f(x)$ by keeping the terms up to $x^2$. Once $\left\langle
x^2\right\rangle M_n^B(x)$ is computed, then the first and the
second moments are easily extracted. In this notation, we have
\[
\left\langle x^2\right\rangle e^{rx}=1+rx+r^2x^2/2.
\]
Moreover,
\begin{align*}
\left\langle x^2\right\rangle\sum_{r=0}^{2j-1}e^{rx}%
&=2j+{2j\choose2}x+\frac{j(2j-1)(4j-1)}{6}x^2,\\[5pt]
\left\langle x^2\right\rangle\prod_{j=k+1}^n\sum_{r=0}^{2j-1}e^{rx}%
&=\left\langle x^2\right\rangle\prod_{j=k+1}^n%
\left(2j+{2j\choose2}x+\frac{j(2j-1)(4j-1)}{6}x^2\right)\\[5pt]
&=\frac{(2n)!!}{(2k)!!}\cdot\left(%
1+\frac{n^2-k^2}{2}x+c_1x^2\right),
\end{align*}
where
\begin{align*}
c_1%
&=\sum_{j=k+1}^n\frac{j(2j-1)(4j-1)}{6}\frac{1}{2j}%
+\sum_{k+1\le i<j\le n}{2i\choose2}{2j\choose2}%
\frac{1}{2i}\frac{1}{2j}\\[5pt]
&=\frac{(n-k)(9n^3+4n^2+9kn^2+6n-9k^2n+4kn-1+6k-9k^3+4k^2)}{72}.
\end{align*}
By the definition \eqref{MGF_B}, we find
\[
M_n^B(x)%
=\frac{1}{D_n^B}\sum_{k=0}^n(-1)^ke^{k(k-1)x}%
\frac{[2n]_{e^x}!!}{[2k]_{e^x}!!}%
=\frac{1}{D_n^B}\sum_{k=0}^n(-1)^ke^{k(k-1)x}%
\prod_{j=k+1}^n\sum_{r=0}^{2j-1}e^{rx}.
\]
It follows that
\begin{align*}
\left\langle x^2\right\rangle M_n^B(x)%
&=\frac{\left\langle x^2\right\rangle}{D_n^B}%
\sum_{k=0}^n(-1)^k\left(1+k(k-1)x+\frac{k^2(k-1)^2}{2}x^2\right)%
\frac{(2n)!!}{(2k)!!}\left(1+\frac{n^2-k^2}{2}x+c_1x^2\right)\\[5pt]
&=\frac{1}{D_n^B}%
\sum_{k=0}^n(-1)^k\frac{(2n)!!}{(2k)!!}%
\left(1+\left(\frac{n^2-k^2}{2}+k(k-1)\right)x+c_2x^2\right),
\end{align*}
where
\[
c_2=c_1+k(k-1)\frac{n^2-k^2}{2}+\frac{k^2(k-1)^2}{2}.
\]
Let $(k)_i=k(k-1)\cdots(k-i+1)$ be the lower factorial. We get
\[
c_2=\frac{9(k)_4+14(k)_3+(18n^2-27)(k)_2-18n^2k+(9n^4+4n^3+6n^2-n)}{72}.
\]
Combining \eqref{D_n_B1}, \eqref{D_n_B3} and \eqref{E_B_tmp}, we
find
\begin{align*}
E_n^B&=[x]M_n^B(x)%
=\frac{1}{D_n^B}\sum_{k=0}^n(-1)^k\frac{(2n)!!}{(2k)!!}%
\left(\frac{n^2-k^2}{2}+k(k-1)\right)\\[5pt]
&=\frac{n^2}{2}+\frac{n}{4}%
+\left(-\frac{n^2}{2}+\frac{3n}{4}\right)\frac{D_{n-1}^B}{D_n^B}.
\end{align*}
Now, the variance of $\fmaj$ equals
\begin{align*}
\frac{72D_n^B}{(2n)!!}[x^2]M_n^B(x)%
&=72\sum_{k=0}^n(-1)^k\frac{c_2}{(2k)!!}\\[5pt]
&=\sum_{k=0}^n(-1)^k%
\left[\frac{9}{2^4(2k-8)!!}%
+\frac{14}{2^3(2k-6)!!}%
+\frac{18n^2-27}{2^2(2k-4)!!}\right.\\[5pt]
&\left.\quad-\frac{18n^2}{2(2k-2)!!}%
+\frac{9n^4+4n^3+6n^2-n}{(2k)!!}\right].
\end{align*}
It can be deduced that
\[
[x^2]M_n^B(x)
=\frac{n(72n^3+140n^2-22n-101)}{576}%
-\frac{n(216n^3-356n^2-186n+127)}{576}\frac{D_{n-1}^B}{D_n^B}.
\]

\begin{thm}\label{EV_B}
The expectation $E_n^B$ and variance $V_n^B$ of $\fmaj$ given by
\[
E_n^B%
=\frac{n^2}{2}+\frac{n}{4}%
+\left(-\frac{n^2}{2}+\frac{3n}{4}\right)\frac{D_{n-1}^B}{D_n^B},
\]
and
\[
V_n^B%
=\frac{n\left(68n^2-40n-101\right)}{288}%
-\frac{n\left(72n^3-212n^2-78n+127\right)}{288}\frac{D_{n-1}^B}{D_n^B}%
-\frac{n^2(2n-3)^2}{16}\left(\frac{D_{n-1}^B}{D_n^B}\right)^2.
\]
\end{thm}

In view of \eqref{D_n_B4} and \eqref{Ratio_D_B}, we obtain the
following estimates.

\begin{cor}\label{EV_Est_B}
We have the following asymptotic estimates:
\[
E_n^B=\frac{n^2}{2}+\frac{3}{8}+o(1),\quad\quad
V_n^B=\frac{n^3}{9}+\frac{n^2}{6}-\frac{n}{36}-\frac{13}{36}+o(1).
\]
\end{cor}

\subsection{The limiting distribution}

We aim to show that the limiting distribution of $\fmaj$ is normal.
The following formula is analogous to Lemma \ref{lem_Tannery_A}.

\begin{lem}\label{lem_Tannery_B}
For any real $x$ satisfying $|x|\le1$ and bounded $|t|< M$,
\[
\lim_{n\to\infty}%
\sum_{k=0}^n\frac{x^k}{[2k]_{e^{-t/\sigma_n^B}}!!\,e^{kt/\sigma_n^B}}%
=e^{x/2}.
\]
\end{lem}

\pf By virtue of Tannery's theorem, if suffices to find an upper
bound $M_k$ for
\[
|v_k(n)|%
=\frac{|x^k|}{[2k]_{e^{-t/\sigma_n^B}}!!\,e^{kt/\sigma_n^B}}%
\]
such that $M_k$ is independent of $n$ and $\sum_{k=0}^{\infty}M_k$
converges.

If $t\le0$, Corollary \ref{EV_Est_B} implies that $\sigma_n^B$ has a
positive lower bound as $n$ runs over all positive integers and so
does $e^{t/\sigma_n^B}$. Suppose that $e^{t/\sigma_n^B}\ge
c_1\in(0,1]$ for all $|t|\le M$, where $c_1$ is independent of $t$.
Then for any $k\ge0$,
\begin{align*}
\frac{|x^k|}{[2k]_{e^{-t/\sigma_n^B}}!!\,e^{kt/\sigma_n^B}}%
&=\prod_{j=1}^k\frac{|x|}%
{\left(1+e^{-t/\sigma_n^B}+e^{-2t/\sigma_n^B}%
+\cdots+e^{-(2j-1)t/\sigma_n^B}\right)e^{t/\sigma_n^B}}\\[5pt]
&=\prod_{j=1}^k\frac{|x|}%
{\left(e^{t/\sigma_n^B}+1+e^{-t/\sigma_n^B}%
+\cdots+e^{-(2j-2)t/\sigma_n^B}\right)}\\[5pt]
&\le(1+c_1)^{-k}.
\end{align*}
Clearly, $\sum_{k\ge0}(1+c_1)^{-k}$ is convergent.

We now assume that $t\ge0$. Suppose $e^{-t/\sigma_n}\ge c_2\in(0,1]$
where $c_2$ is independent of $t$. Then for any $k\ge1$,
\begin{align*}
\frac{|x^k|}{[2k]_{e^{-t/\sigma_n^B}}!!\,e^{kt/\sigma_n^B}}%
&=\prod_{j=1}^k\frac{|x|}%
{\left(1+e^{-t/\sigma_n^B}+\cdots+e^{-(2j-1)t/\sigma_n^B}\right)%
e^{kt/\sigma_n^B}}\\[5pt]
&\le\prod_{j=1}^k\frac{1}{1+c_2+\cdots+c_2^{2j-1}}%
\le(1+c_2)^{-k}.
\end{align*}
Similarly, $\sum_{k\ge0}(1+c_2)^{-k}$ is convergent. \qed

The following formula, which is similar to Lemma \ref{lem_i>=2_A},
will be crucial in the proof of main theorem of this section.

\begin{lem}\label{lem_i>=2_B}
For any bounded $|t|<M$,
\begin{equation}\label{series}
\lim_{n\to\infty}\sum_{i=2}^{\infty}%
\frac{B_{2i}\,t^{2i}}{(2i)\,(2i)!\,\left(\sigma_n^B\right)^{2i}}%
\sum_{j=1}^n\left((2j)^{2i}-1\right)=0,
\end{equation}
where $B_{2i}$ are the Bernoulli numbers.
\end{lem}

\pf Let $\alpha$, $\beta$ and $\gamma$ be three constants such that
$\alpha>2$, $\beta>9$, and $0<\gamma<1/2$. Let $N$ be a fixed
integer satisfying the following three conditions:
\begin{itemize}
\item%
$2n+2<\alpha n$ for any $n>N$;
\item%
$\left(\sigma_n^B\right)^2-n^3/\beta>0$ for any $n>N$;
\item%
$2\pi N^{\gamma/2}>M\alpha\sqrt{\beta}$.
\end{itemize}
The existence of such $N$ is evident. Let $i\ge2$ and $n>N$. We will
show that the series in \eqref{series} is convergent to zero
absolutely. It is easy to derive the following upper bound:
\[
\sum_{j=1}^{n}\left((2j)^{2i}-1\right)%
=2^{2i}\sum_{j=1}^nj^{2i}-n <2^{2i}\int_1^{n+1}t^{2i}dt\]
\[
<\frac{2^{2i}(n+1)^{2i+1}}{2i+1}%
<\frac{(2n+2)^{2i+1}}{5}%
<\frac{(\alpha n)^{2i+1}}{5}.
\]
The rest of the proof is similar to that of Lemma \ref{lem_i>=2_A}.
\qed

\begin{lem}\label{lem_MGF_B}
The following relation holds:
\[
M_n^B(x)=\frac{(2n)!!}{D_n^B}\,\exp\left(\frac{xn^2}{2}%
+\sum_{i=1}^\infty\frac{B_{2i}\,x^{2i}}{(2i)\,(2i)!}%
\sum_{j=1}^n\left((2j)^{2i}-1\right)\right)%
\sum_{k=0}^n\frac{(-1)^k}{[2k]_{e^{-x}}!!\,e^{kx}}.
\]
\end{lem}

\pf From \eqref{MGF_B} and \eqref{d_n(q)_B}, we have
\begin{equation}\label{eq5}
M_n^B(x)=d_n^B\left(e^x\right)/D_n^B
=\frac{1}{D_n^B}\,[2n]_{e^x}!!%
\sum_{k=0}^n\frac{(-1)^k\,e^{xk(k-1)}}{[2k]_{e^x}!!}.
\end{equation}
Moreover,
\begin{align}
[2n]_{e^x}!!%
&=\prod_{j=1}^n[2j]_{e^x}%
=\prod_{j=1}^n(2j)\cdotp\exp\left(\frac{x(2j-1)}{2}%
+\sum_{i=1}^\infty%
\frac{B_{2i}\,x^{2i}\left((2j)^{2i}-1\right)}{(2i)\,(2i)!}%
\right)\notag\\[5pt]
&=(2n)!!\cdotp%
\exp\left(\frac{xn^2}{2}+%
\sum_{i=1}^\infty\frac{B_{2i}\,x^{2i}}{(2i)\,(2i)!}%
\sum_{j=1}^n\left((2j)^{2i}-1\right)\right),\label{eq6}
\end{align}
and
\begin{equation}\label{eq7}
\frac{e^{xk^2}}{[2k]_{e^x}!!}%
=\prod_{j=1}^k\frac{e^{(2j-1)x}}{[2j]_{e^x}}%
=\prod_{j=1}^k\frac{e^{2jx}}{e^x}\frac{1-e^x}{1-e^{2jx}}%
=\prod_{j=1}^k\frac{1-e^{-x}}{1-e^{-2jx}}%
=\frac{1}{[2k]_{e^{-x}}!!}.
\end{equation}
Substituting \eqref{eq6} and \eqref{eq7} into \eqref{eq5}, we deduce
the desired relation. \qed

\begin{thm}\label{Main_A}
Let $\fmaj$ be the flag major index of a random $B_n$-derangement.
Then the distribution of the random variable
\[
\xi_n^B={\fmaj-E_n^B \over \sigma_n^B}
\]
converges to the standard normal distribution as $n\to\infty$.
\end{thm}

 \pf By Curtiss's theorem and \eqref{MGF_Standard_B}, normality of
the distribution of $\fmaj$ follows from the identity
\begin{equation}\label{Proof_Main_B}
\lim_{n\to\infty}e^{-t\,E_n^B/\sigma_n^B}M_n^B\left(t/\sigma_n^B\right)%
=e^{t^2/2}.
\end{equation}
By Lemma \ref{lem_MGF_B}, the left hand side of \eqref{Proof_Main_B}
can be expressed as the limit of the following expression:
\[
\frac{(2n)!!}{D_n^B}%
\exp\left(-\frac{tE_n^B}{\sigma_n^B}+\frac{tn^2}{2\sigma_n^B}%
+\sum_{i=1}^\infty\frac{B_{2i}\,t^{2i}}%
{(2i)\,(2i)!\left(\sigma_n^B\right)^{2i}}%
\sum_{j=1}^n\left((2j)^{2i}-1\right)\right)%
\sum_{k=0}^n\frac{(-1)^k}{[2k]_{e^{-t/\sigma_n^B}}!!\,e^{kt/\sigma_n^B}}.
\]
First, the estimate \eqref{D_n_B4} implies that
\begin{equation}\label{eq1_B}
\lim_{n\to\infty}(2n)!!/D_n^B=\sqrt{e}.
\end{equation}
By Corollary \ref{EV_Est_B}, for bounded $t$ we have
\begin{equation}\label{eq2_B}
\lim_{n\to\infty}%
\left(-\frac{tE_n^B}{\sigma_n^B}+\frac{tn^2}{2\sigma_n^B}\right)%
=\lim_{n\to\infty}%
\frac{t}{\sigma_n^B}\left(\frac{n^2}{2}-E_n^B\right)%
=0.
\end{equation}
It is easy to check that
\[
\lim_{n\to\infty}%
\frac{1}%
{\left(\sigma_n^B\right)^{2}}%
\sum_{j=1}^n\left((2j)^{2}-1\right)=12.
\]
Based on Lemma \ref{lem_i>=2_B} and the fact that $B_2=1/6$, we see
that
\begin{equation}\label{eq3_B}
\lim_{n\to\infty}%
\sum_{i=1}^\infty\frac{B_{2i}\,t^{2i}}%
{(2i)\,(2i)!\left(\sigma_n^B\right)^{2i}}%
\sum_{j=1}^n\left((2j)^{2i}-1\right)%
=\lim_{n\to\infty}%
\frac{B_2t^2}{2\cdot2!\,\left(\sigma_n^B\right)^{2}}%
\sum_{j=1}^n((2j)^2-1)%
=\frac{t^2}{2}.
\end{equation}
Finally, taking $x=-1$ in Lemma \ref{lem_Tannery_B}, we get
\begin{equation}\label{eq4_B}
\lim_{n\to\infty}%
\sum_{k=0}^n\frac{(-1)^k}{[2k]_{e^{-t/\sigma_n^B}}!!\,e^{kt/\sigma_n^B}}%
=e^{-1/2}.
\end{equation}
Combining \eqref{eq1_B}, \eqref{eq2_B}, \eqref{eq3_B} and
\eqref{eq4_B}, we obtain \eqref{Proof_Main_B}. This completes the
proof. \qed

\noindent{\bf Acknowledgments.} This work was supported by the 973
Project, the PCSIRT Project of the Ministry of Education, the
Ministry of Science and Technology, and the National Science
Foundation of China.


\end{document}